\documentclass[12pt]{article}

 \usepackage{amsfonts,amssymb,graphicx}

\newtheorem{prop}{Proposition}
\newtheorem{rema}{Remark}
\newtheorem{defi}{Definition}
\newtheorem{lemm}{Lemma}

\newcommand{\C}[1][]{\ensuremath{{\mathbb{C}^{#1}} }}
\newcommand{\R}[1][]{\ensuremath{{\mathbb{R}^{#1}} }}

\renewcommand{\S}[1][]{\ensuremath{{\mathbb{S}^{#1}} }}

\def\Re{ \mathrm{Re}\, }
\def\Im{ \mathrm{Im}\, }

\newcommand{\ga}{\gamma}
\newcommand{\al}{\alpha}
\newcommand{\be}{\beta}
\newcommand{\pa}{\partial}

\def\mail#1#2{{\tt #1}@{\tt #2}}
\date{}
\title{Special Lagrangian submanifolds in the complex sphere}
\author{Henri Anciaux}

\begin{document}

\maketitle

\vspace{-2.4em}
\centerline{ \small email: \mail{henri}{impa.br}}
\vspace{3em}
\centerline{\textbf {Abstract}}

\medskip

{ \small We construct a family of Lagrangian submanifolds in the complex 
shere whith a $SO(n)$-invariance property. Among them we find those which
are special Lagrangian with respect with the Calabi-Yau structure defined by
 the
Stenzel metric}

\centerline{\small \em MSC: 53D12, 53C42, 49Q05 \em}

\bigskip

\section*{Introduction}
Special Lagrangian submanifolds may be defined as those
submanifolds which are both Lagrangian (an order 1 condition) and
minimal (an order 2 condition). Alternatively, they are
characterised as those submanifolds which are calibrated by a
certain $n$-form (\em cf \em [HL]), so they have the remarkable
property of being area minimizing.
 Their study have received many
attention recently since connections with string theory have been
discovered. More particularly, understanding fibrations of special
Lagrangian (possibly with singularities) in Calabi-Yau manifolds
of (complex) dimension 3 is crucial for mirror symmetry (\em cf
\em [SYZ],[Jo]). Since the pioneering work of Harvey and Lawson
[HL], where this notion  
 were introduced, several authors
([Ha],[Jo],[CU2]) have discovered many families of special
Lagrangian submanifolds in the complex Euclidean space, however
very few examples (\em cf \em [Br]) of such submanifolds are known
in other Calabi-Yau manifolds, where this notion is naturally
extended. Maybe the main reason is that such manifolds are
somewhat rare, in particular the existence of compact, Calabi-Yau
manifolds is a hard result of S.-T. Yau [Y] involving non explicit
solutions to some PDE.   
However there exists some intermediary
examples of non-flat, non-compact Calabi-Yau manifolds, maybe the
simplest of them being the complex sphere.

 In this paper we describe a family of Lagrangian submanifolds
 of a complex variety of $\C^{n+1}$ with a $SO(n)$-invariance.
In the case of the complex sphere, we obtain among them new
examples of
 special Lagrangian
 submanifolds, which are a kind of
 generalization of the Lagrangian catenoid discovered in [HaLa]
 (\em cf \em also [CU1]).

 In [CMU], minimal Lagrangian submanifolds have been studied in
the complex hyperbolic and complex projective spaces respectively.
This is somewhat close to our situation for the following reason:
these spaces are not Calabi-Yau manifolds, and so these
submanifolds are not minimizing \em a priori
\em; however Y.-G. Oh
has proved in [Oh] that a minimal Lagrangian submanifold of the
complex hyperbolic space is stable and without nullity, a property
necessary for being a minimizer.

The paper is organized as follows:

In Section 1 we give a definition of Calabi-Yau manifolds, we
state some basic facts on them and describe the Calabi-Yau
structure yielded on the complex sphere by the Stenzel metric. In
Section 2 we describe a general class of Lagrangian submanifolds
whith some $SO(n)$-invariance immersed in some complex submanifold
of $\C^{n+1}.$ In the particular case of the complex sphere, we
then check that they are also Lagrangian for the Calabi-Yau
structure. In the last section we find in this class the special
Lagrangian and describe them.

The nice argument about
$SO(n+1)$-invariance at the end of next section is due to
Dominic Joyce and Mu-Tao Wang
suggested to me the Lagrangian Ansatz; it is my pleasure to
thank them. The author, as a member of EDGE, 
Research Training Network HPRN-CT-2000-00101,
has been supported by the European Human Potential Programme.

\section{The Calabi-Yau structure on the complex sphere}
We start with the following definitions --- which are not
the standard ones --- proposed in [Jo]:

\begin{defi}
 Let $n \geq 2.$ An almost Calabi-Yau manifold is a quadruple
 $(X,J,\omega,\Omega)$ such that $(X,J)$ is a
 $n$-dimensional complex manifold, $\omega$ the symplectic form of a
 K\"{a}hler metric $g$ on $X$ and $\Omega$
 a non-vanishing holomorphic
$(n,0)$-form.

$(X,J,\omega,\Omega)$ is called a Calabi-Yau manifold if in
addition
$$ \frac{\omega^n}{n!}= (-1)^{n(n-1)/2} \Omega \wedge \bar{\Omega}, \leqno{(1)}$$
 where $\omega$ is the symplectic $2$-form corresponding to $g.$

\end{defi}

An important property of a Calabi-Yau manifold is that it is
Ricci-flat, \em i.e. \em the Ricci curvature vanishes. Conversely,
on a simply connected open subset a Ricci-flat K\"{a}hler metric
yields a Calabi-Yau structure. The most simple example of
Calabi-Yau manifold is $\C^n$ equipped with its standard
structures, that is
$$ \omega_{0}= i \sum_{j} dz_{j} \wedge d\bar{z}_{j}$$
and
$$ \Omega_0= dz_1 \wedge \ldots \wedge dz_n.$$

\begin{defi} A submanifold $L$ of real dimension $n$ of $X$
is said to be special Lagrangian if $\omega|_{L}\equiv \Im \Omega
|_{L} \equiv 0.$
\end{defi}

The first condition says that $L$ is Lagrangian and the second one
that $L$ is calibrated with respect to $\Re \Omega.$ This implies
that a special Lagrangian is always necessarily minimizing (\em cf
\em [HL]).

\bigskip

We now describe a class of almost Calabi-Yau manifolds. Let $P$ be
a complex polynomial. Then the set
$$ Q:=\{ (z_{0}, \ldots ,z_{n}), P(z_{0})= \sum_{j=1}^{n} z_{j}^{2} \} $$
 is a complex submanifold of $\C^{n+1}$ which is smooth except if
 $P$ admits double poles. For example, if $P(z)=z^2,$ $Q$ is
 made of two copies of $\C^n$ singularly connected at $0.$
$Q$ inherits from $\C^{n+1}$ a K\"{a}hler
 structure.

 On $Q$ we define the following holomorphic
 $(n,0)$-form $\Omega,$ defined by
$$ \Omega \wedge d(-P(z_0)+ z_1^2 + \ldots + z_n^2)= dz_0 \wedge
\ldots \wedge dz_n$$
 We may give an explicit form on the set $\{P(z_0) \neq 0  \}:$

 $$ \Omega= \frac{1}{\dot{P}(z_0)} dz_1 \wedge \ldots \wedge dz_n $$

In the remainder of the section we shall see that in the
particular case of the complex sphere, that is when $P(z)=1-z^2,$
$Q$ is a Calabi-Yau manifold.
 This is a consequence of
the following result due to Stenzel (\em cf \em [St]):

\begin{prop}
 On the complex sphere $\{\sum_{j=0}^{n} z_{j}^{2}=1 \} ,$
 there exists a Ricci-flat metric, whose corresponding symplectic
 form is
$$ \omega_{St}=i \pa \bar{\pa} u(r^{2})=
  i \sum_{j,k} \frac{\partial^2}{\partial z_j \partial \bar{z}_k} u(r^2) dz_j \wedge d\bar{z}_k$$
  where $r^{2}= \sum_{j=0}^{n} z_j \bar{z}_j$ and $u$ is some smooth
  real function.
\end{prop}

  In dimension 2, we have an explicit form: $u(r^{2})=\sqrt{1+r^{2}},$
and in higher dimension, $u$ is defined as the solution of some
differential equation. As the complex sphere is simply connected,
this result shows that it is Calabi-Yau.

It remains to show that the holomorphic $n$-form $\Omega$ that we
have introduced is the one of the Calabi-Yau structure, up to a
multiplicative constant. So see this, let $f$ be the ratio of
these two holomorphic forms. We observe that $\Omega$ is invariant
under the natural action of $SO(n+1)$ defined by $z=x+iy \mapsto
Az=Ax+iAy, \forall A \in SO(n+1).$ But this is also the case for
the Stenzel metric, and thus for the corresponding holomorphic
$n$-form. This implies that $f$ is constant on the orbits of the
action of $SO(n+1).$ But generic orbits have codimension $1$ in
the complex sphere, so the image of $f$ in $\C$ has at most real
dimension 1 locally, which forces $f$ to be constant, as it is
holomorphic. Thus the two forms are proportional.

\section{A Lagrangian Ansatz}

Let $\gamma$ be a curve into $\C.$
 Then
   $$ \begin{array}{lccc} X : & I \times \S^{n-1}  & \to & Q  \\
&  (s , x) & \mapsto &
 (\ga(s), \sqrt{P(\ga(s))} x_1, \ldots , \sqrt{P(\ga(s))} x_n ),
\end{array}$$
 where $(x_1, \ldots, x_n)=x,$
 is a an Lagrangian immersion for $\omega_{0}.$
In spite of the indetermination of the complex square root, the
immersion is well defined because of the invariance of $\S^{n-1}$
by $x \mapsto -x.$ However it becomes singular when $\ga$ is a
zero of $P.$

\bigskip

From now on we consider only the case $P(\ga)=1- \ga^2,$ that is
the complex sphere equipped with the Stenzel form

Due to the larger symmetry of the complex sphere, we may slightly
generalise the previous ansatz: we embed $\S^{n-1}$
 into $\S^{n}$ in the following way: $\S^{n-1}= \S^{n} \cap \{ x_0=0\}.$
 For some $p \in \S^n,$ we shall note $A_p$ any element of $SO(n+1)$ such that
 $A_p(\{ x_0=0\})=p^{\perp}.$
Let $\gamma$ be a curve into $\C^{*}.$
 Then we define
   $$ \begin{array}{lccc} X_{p,\ga} : & I \times \S^{n-1}  & \to & Q  \\
&  (s , x) & \mapsto & \ga(s) p + \sqrt{1-\ga(s)^2} A_p x,
\end{array}$$
 and it is clear that the image does not depend on the choice of
 $A_p.$

   \begin{lemm} $X$ is also Lagrangian for the Stenzel metric. \end{lemm}

 \noindent \textbf{Proof.}
 Let $(z_0, \ldots , z_n)$ be coordinates on $\C^n$ such that
 $p=(1,0,\ldots,0).$ As long as $\ga $ does not vanish
 we can use $(z_1, \ldots,z_n)$ as coordinates
 on a neighbourhood of the image of $X$ in $Q.$
   In particular, we have
  $$ \frac{\partial z_0}{\partial z_j}=-\frac{z_j}{z_0}.$$

 We compute that
  $$ \frac{\partial}{\partial z_j} r^2= \bar{z_j} -
  \frac{\bar{z}_0}{z_0} z_j $$
  $$  \frac{\partial^2}{\partial z_j \partial \bar{z}_k} r^2=
       \delta_{jk}+\frac{z_j \bar{z}_k}{|z_0|^2}$$
  So we deduce that
  $$\omega_{St}=i \sum_{j,k}^{n} a_{jk} dz_j \wedge d\bar{z}_k$$
  with
  $$a_{jk}=\left(\delta_{jk}+\frac{z_j \bar{z}_k}{|z_0|^2}\right) u' +
   2 \Re \left(\bar{z}_j  z_k - \frac{\bar{z}_0}{z_0} z_j z_k \right) u'' $$
Hence we can decompose $\omega_{St}$ as
  $$\omega_{St} = u'  \omega_0 + \omega_1,$$
where
 $$\omega_1:=  i\sum_{j,k}^{n} \left(\frac{z_j \bar{z}_k}{|z_0|^2} u' +
   2 \Re \left(\bar{z}_j  z_k - \frac{\bar{z}_0}{z_0} z_j z_k \right) u''\right) dz_j \wedge d\bar{z}_k$$
As $X$ is Lagrangian for $\omega_0,$ it remains to show that it is
also Lagrangian for $\omega_1.$ This will be a consequence of the
fact that the two following $2$-forms vanish on $X:$
 $$ \sum_{j,k}^{n} z_j \bar{z}_k  dz_j \wedge d\bar{z}_k $$
 $$ \sum_{j,k}^{n} z_j z_k  dz_j \wedge d\bar{z}_k $$
We shall note $T_\al = (X_{\al}^1, \ldots , X_{\al}^n),$
 $ 1 \leq \al \leq n-1,$
 a basis of tangent vectors to $\S^{n-1}$ at
 $x=(x_1, \ldots , x_n).$

 This yields a basis of the tangent space at $X(s,x)$:
 $$ X_s=(\dot{\ga},  \frac{\dot{\ga}\dot{P}(\ga)}{2 \sqrt{P(\ga)}} x_1, \ldots ,
                  \frac{\dot{\ga}\dot{P}(\ga)}{2 \sqrt{P(\ga)}} x_n)$$
  $$ X_{\al}=(0, \sqrt{P(\ga)} X_{\al}^{1}, \ldots ,
                       \sqrt{P(\ga)} X_{\al}^{n}).$$

\bigskip

Now we have:
$$  \left(\sum_{j,k} \bar{z}_j z_k dz_j \wedge d\bar{z}_k \right)(X_s,X_{\al})
 = \sum_{j,k} |P(\ga)| x_j x_k
     \left| \begin{array}{cc}
       \frac{\dot{\ga}\dot{P}(\ga)}{2 \sqrt{P(\ga)}} x_j &
       \sqrt{P(\ga)}X_j^{\al}\\ & \\
     \overline{ \frac{\dot{\ga}\dot{P}(\ga)}{2 \sqrt{P(\ga)}} } x_k &
       \overline{\sqrt{P(\ga)}} X_k^{\al}
         \end{array}  \right| $$
 $$= \sum_{j,k} |P(\ga)| x_j x_k
   \left(  \frac{\dot{\ga}\dot{P}(\ga)}{2 \sqrt{P(\ga)}} \overline{\sqrt{P(\ga)}}
       x_j X_\al^k
   - \overline{\frac{\dot{\ga}\dot{P}(\ga)}{2 \sqrt{P(\ga)}}} \sqrt{P(\ga)}
      x_k X_\al^j
       \right) $$
  $$=   |P(\ga)|
     \frac{\dot{\ga} \dot{P}(\ga)}{2 \sqrt{P(\ga)}} \overline{\sqrt{P(\ga)}}
    \sum_j x_j^{2} \sum_k x_k X_\al^k
   -  |P(\ga)|
   \overline{\frac{\dot{\ga}\dot{P}(\ga)}{2 \sqrt{P(\ga)}}} \sqrt{P(\ga)}
    \sum_k x_k^{2} \sum_j x_j X_\al^j $$
The latter vanish because tangent vectors to the real sphere at
some point are orthogonal to this point.

 \bigskip
$$  \left(\sum_{j,k} z_j \bar{z}_k  dz_j \wedge d\bar{z}_k \right)(X_{\al},X_{\be})
=  \sum_{j,k}  |P(\ga)| x_j x_k
  \left| \begin{array}{cc}
       \sqrt{P(\ga)} X_j^{\al} &
       \sqrt{P(\ga)}X_j^{\be}\\ & \\
     \overline{ \sqrt{P(\ga)}} X_k^{\al} &
       \overline{\sqrt{P(\ga)}} X_k^{\be}
         \end{array}  \right| $$
   $$=\sum_{j,k} |P(\ga)| x_j x_k \left(
    P(\ga)X_j^{\al} X_k^{\be} -P(\ga) X_k^{\al} X_j^{\be}\right)=0$$

The computations for showing that also $\sum_{j,k} z_j \bar{z}_k
dz_j \wedge dz_k$ vanishes on $X$ are analogous.

\section{Special Lagrangian submanifolds}
 Into the class of Lagrangian defined in the previous section,
and in the two cases where we have a Calabi-Yau structure,
 we
 shall look for those calibrated by $\Re \Omega,$ that is those
 on which $\Im \Omega$ vanishes.

\subsection{Flat case of $P(z)=z^2$}
 The holomorphic $(n,0)$-form we use is the standard one
$$\Omega_0= dz_1 \wedge \ldots \wedge dz_n$$

The equation is $\Im( \dot{\ga} \ga^{n-1}) =0,$ which is easily
integrated as $\Im(\ga^n)=C,$ where $C$ is some real constant. It
turns that we have two types of solutions:

\begin{itemize}

\item If $C=0,$ $\ga$ is the union of $n$ lines passing trough the
origin, with angles $\pi / n [2 \pi]$ and the corresponding
special Lagrangian are just Lagrangian $n$-spaces,

\item If $C \neq 0,$ the curve is asymptotic to two of the
lines described above. All the curves (and so the corresponding
Lagrangian as well) are congruent. We recognize the Lagrangian
catenoid which was first identified in [HL] and characterised in
[CU1] as the only special Lagrangian of $\C^n$ which is foliated
by $(n-1)$-dimensional round spheres.

\end{itemize}

\subsection{Case of the complex sphere $P(z)=1-z^2$}
Before starting the computations, we observe that we already know
a special Lagrangian in the complex sphere: if a Calabi-Yau has a
real structure which is in some sense compatible (\em cf \em [Au],
pp. 62--64), then the set of real points, if non empty, is a
special Lagrangian submanifold. In our case, real points
constitute an real sphere embedded in $Q,$ so this will be a
trivial solution to our problem.

\bigskip

 We now compute the holomorphic form on $X:$

$$\Omega(X_s,X_1, \ldots , X_{n-1})=\frac{1}{\dot{P}(\ga)} \left|
   \begin{array}{cccc}
           \frac{\dot{\ga}\dot{P}(\ga)}{2\sqrt{P(\ga)}}x_1
                                            & \sqrt{P}X_1^1 & \cdots & \sqrt{P} X_{n-1}^1 \\
    \vdots & \vdots & \ddots & \vdots \\
   \frac{\dot{\ga}\dot{P}(\ga)}{2 \sqrt{P(\ga)}} x_n
                                 &  \sqrt{P}X_1^n & \cdots &  \sqrt{P}X_{n-1}^n  \end{array}
                       \right|$$
            $$ =  \frac{\dot{\ga}}{2\sqrt{P(\ga)}}
            \sqrt{P(\ga)}^{n-1}  \left|
   \begin{array}{cccc}
           x_1                & X_1^1 & \cdots &  X_{n-1}^1 \\
    \vdots & \vdots & \ddots & \vdots \\
    x_n           &   X_1^n & \cdots &  X_{n-1}^n  \end{array}  \right|$$
                        $$ =  \dot{\ga} \sqrt{P(\ga)}^{n-2}  C_{\R},$$
 where $C_{\R}$ is some real constant.

 Thus the differential equation for $\ga$ is
 $$ \Im \left( \dot{\ga} \sqrt{1-\ga^2}^{n-2} \right)= 0.$$

This equation has two singular points $\pm 1$ and is regular
elsewhere. Moreover, there is always a simple solution to this
equation, the real
 segment $[-1,1].$ This corresponds to the
standard embedding of the real sphere. We shall analyse separately
the even and odd cases.

\subsubsection{Even case}
 In this case, the equation is polynomial and it is easy to find a first integral:
 we integrate $\Im \left( \dot{\ga}
(1-\ga^2)^{n/2-1} \right) =0$ and we get
 $ \Im(Q(\ga))= C,$ where $Q$ is some polynomial such that
 $Q'(z)=(1-z^2)^{n/2-1}.$ The degree of $Q$ is $n-1$ and the the integral
 curves
are algebraic curves of degree $n-1$ as well.

In order to have a better description of the solutions, we perform
an asymptotic analysis of the equation when $|\ga| \to \infty$ and
when $\ga \sim \pm 1.$

If $|\ga| > > 1,$ we have
   $$\Im \left( \dot{\ga} (1-\ga^2)^{n/2-1} \right) \sim -\Im(\dot{\ga} \ga^{n-2}), $$
so asymptotically, the phase portrait of the equation looks like
the one of the flat case (however not of the same dimension). In
particular all integral curves are asymptotic to the half lines
$\{ \arg(z)=  k \pi / (n-1) \}$ , $1 \leq k \leq 2n-2.$

Next we write $\ga= 1 + y$ and we find that for $y$ close to $0,$
$$\Im \left( \dot{\ga} (1-\ga^2)^{n/2-1} \right) \sim
             \Im \left( \dot{y} (- 2y)^{n/2-1}\right)$$
So close to the singular point $1,$ the integral lines look like
the level sets of
$$ \Im(-1/n (-2y)^{n/2})=C.$$
 In particular, the level set of level
$0$ is the union of $n$ branches asymptotic in $1$ to the
half-lines $\{ \arg(z)= 2k \pi /n+ \pi \}$ , $1 \leq k \leq n,$
one of them being of course the real segment $[-1,1].$

At the point $-1$ the situation is exactly symmetric with respect
to the reflexion $ z \to -\bar{z}.$

\subsubsection{Odd case}

When $n$ is odd, there is also a first integral which takes the
following form
 $$ \Im \left(\ga \sqrt{1-\ga^2} R(\ga) + A \arcsin(\ga) \right),$$
where $R$ is a polynomial of degree $n-3$ and $A$ a real constant.
This quantity is not well-defined on the whole complex plane,
however the level set lines can be defined locally.

If $|\ga|  >> 1,$ we have
   $$\Im \left( \dot{\ga} \sqrt{1-\ga^2}^{n-2} \right)
    =\Im \left( \dot{\ga} \sqrt{1+(i\ga)^2}^{n-2} \right)
   \sim \Im( \dot{\ga} (i\ga)^{n-2} )
   =\pm \Im(i \dot{\ga} \ga^{n-2} ) $$
Again, the phase portrait looks asymptotically as in the flat
case, but this time up to a rotation of angle $\pm \pi / 2(n-1).$ In
particular, all curves are asymptotic to the half lines
$\{\arg(z)=  k \pi / (n-1)+ \pi / 2(n-1) \}$ , $1 \leq k \leq 2n-2.$

The asymptotic analysis next to the singular points $\pm 1$ is the
same than in the even case.

\subsubsection{Conclusion}
 From these remarks, we can describe the general picture of the
 phase portrait :
\begin{itemize}
\item There is a singular curve made of the real segment $[-1,1]$
and of $2n-2$ branches, one half of them starting from $1$ and the
other half from $-1,$ making between them an angle of $2 \pi/ n$
and going to infinity where they are asymptotic to the half lines
$\{ \arg(z)=  k \pi / (n-1) \}$ , $1 \leq k \leq 2n-2$ when $n$ is
even and to $\{ \arg(z)=  k \pi / (n-1)+ \pi / 2(n-1) \}$ , $1 \leq k
\leq 2n-2$ when $n$ is odd.

\item the other curves are smooth and have two ends asymptotic to
two successive branches described above
\end{itemize}

 For the corresponding special Lagrangian in $Q,$ we deduce the
 following:

\begin{itemize}

\item The standard embedding of the real sphere,
\item $2n-2$ special Lagrangian,
 touching  the real sphere at the north pole $(1,0, \ldots, 0)$
 for one half of them and at the south pole $(-1,0 \ldots, 0)$ for
 the other half.
\item A one-parameter smooth families of special Lagrangian with two
ends asymptotic to two members of the above family.

\end{itemize}

\begin{rema}
In the case of dimension $2,$ the integral lines are simply the
horizontal lines.
\end{rema}

%









  \vspace{2cm}

\noindent Henri Anciaux,

\noindent IMPA,   Estrada Dona Castorina, 110 

\noindent 22460--320 Rio de Janeiro, Brasil
\end{document}